\documentclass[12pt]{amsart}

\newcommand{\szego}{Szeg\"o }

\newcommand{\kahler}{K\"ahler }

\newcommand{\R}{{\mathbb R}}
\newcommand{\C}{{\mathbb C}}

\newcommand{\Z}{{\mathbb Z}}

\newcommand{\Ss}{{\mathbb S}}
\newcommand{\Hh}{{\mathbb H}}

\newcommand{\dbar}{\bar\partial}
\newcommand{\ddbar}{\partial\dbar}

\newcommand{\E}{{\mathbf E}}

\renewcommand{\phi}{\varphi}

\newcommand{\dcal}{\mathcal{D}}

\newcommand{\hcal}{\mathcal{H}}

\newcommand{\ocal}{\mathcal{O}}

\renewcommand{\phi}{\varphi}

\newtheorem{theo}{{\sc Theorem}}[section]
\newtheorem{cor}[theo]{{\sc Corollary}}
\newtheorem{defn}[theo]{{\sc Definition}}

\newtheorem{conj}[theo]{{\sc Conjecture}}
\newtheorem{lem}[theo]{{\sc Lemma}}

\newtheorem{prop}[theo]{{\sc Proposition}}

\title[Complex zeros of real ergodic eigenfunctions]{Complex  zeros of real ergodic eigenfunctions}

\author{Steve Zelditch}
\address{Department of Mathematics, Johns Hopkins University, Baltimore, MD
21218, USA} \email{ zelditch@@math.jhu.edu }

\thanks{Research  partially supported by  NSF grants \#DMS 0302518 and  \# FRG 0354386 .}

\date{August 5, 2005}

\begin{document}

\begin{abstract} We determine the limit distribution (as $\lambda \to \infty$) of complex zeros for
holomorphic continuations $\phi_{\lambda}^{\C}$ to Grauert tubes
of real eigenfunctions of the Laplacian on a real analytic compact
Riemannian manifold $(M, g)$ with ergodic geodesic flow. If
$\{\phi_{j_k} \}$ is an ergodic sequence of eigenfunctions, we
prove the weak limit formula $\frac{1}{\lambda_j}
[Z_{\phi_{j_k}^{\C}}] \to \frac{i}{ \pi} \overline{\partial}
{\partial} |\xi|_g$, where $ [Z_{\phi_{j_k}^{\C}}]$ is the current
of integration over the complex zeros and where
$\overline{\partial}$ is with respect to the adapted complex
structure of Lempert-Sz\"oke and Guillemin-Stenzel.
\end{abstract}

\maketitle

\section{Introduction} A well-known problem in the geometry of
Laplace  eigenfunctions is to determine the asymptotics of the
volume of their nodal hypersurfaces in the limit of large
eigenvalues. At the present time, the best result for compact real
analytic Riemannian manifolds $(M, g)$ of dimension $m$ is the
estimate
\begin{equation} \label{DF} c_1 \lambda \leq {\mathcal
H}^{m-1}(Z_{\phi_{\lambda}}) \leq C_2 \lambda. \end{equation} due
to Donnelly-Fefferman \cite{DF}. Further background and references
can be found in \cite{DF, JL, NPSo}.  In this article we are
concerned with the yet more difficult problem of the  asymptotic
distribution of nodal hypersurfaces, i.e. with integrals of
continuous functions over nodal hypersurfaces (cf.
(\ref{REALINT})). This problem is much too difficult for real
hypersurfaces, but it turns out to simplify quite a bit if we
complexify the problem, i.e. analytically continue the
eigenfunctions into the complexification of $M$.  Our main results
determine the asymptotic distribution of complex nodal
hypersurfaces for eigenfunctions of real analytic Riemannian
manifolds with ergodic geodesic flow, or more generally for any
sequence of quantum ergodic eigenfunctions.

To state our results, we need some notation. Let $(M, g)$ be   a
real analytic Riemannian manifold of dimension $m$, and
 consider an orthonormal basis of real eigenfunctions
$$\Delta_g \phi_{j} = \lambda_j^2 \phi_j, \;\;\;
\langle \phi_j, \phi_k \rangle_{L^2(M, dvol_g)}  = \int_M
\phi_j(x) \phi_k(x) dvol_g(x) =  \delta_{jk}$$ of its  Laplacian
$\Delta_g$. We use the sign convention for which the  Laplacian is
positive; we often write $\Delta_g$ as $\Delta$ when the metric is
understood.
 As reviewed in \S \ref{AC}, a real analytic
manifold possesses a Bruhat-Whitney complexification $M_{\C}$,
that is, a complex manifold in which $M$ embeds as  a totally real
submanifold.  This complex manifold may be identified (by means of
the complexified exponential map) with a ball bundle
$B^*_{\epsilon_0} M$ inside the cotangent bundle $T^*M$, equipped
with a complex structure $J_g$ adapted to the metric in the sense
of Guillemin-Stenzel \cite{GS1, GS2} and Lempert-Sz\"oke
\cite{LS1, LS2}.  We denote by $|\xi|_g^2 = \sum_{ij} g^{ij} \xi_i
\xi_j$ the norm-squared of $\xi \in T^*_x M$  with respect to the
metric, and by $\overline{\partial}$ the Cauchy-Riemann operator
with respect to $J_g$. The maximal ball bundle $B^*_{\epsilon_0}
M$ to which this complex structure extends is known as the Grauert
tube of $(M, g)$. The natural symplectic form $\omega$ and the
complex structure $J_g$ endow $B^*_{\epsilon_0}M$ with the \kahler
metric $\omega_g =  \frac{1}{i} \overline{\partial}
\partial |\xi|_g^2$.

By a theorem due to Boutet de Monvel \cite{Bou} (see also
\cite{GS2, GLS}) the  eigenfunctions possess analytic
continuations $\phi_{\lambda}^{\C}$ to the maximal Grauert tube.
The complex nodal hypersurface of an eigenfunction is   defined by
\begin{equation} Z_{\phi_{\lambda}^{\C}} = \{\zeta \in
B^*_{\epsilon_0} M: \phi_{\lambda}^{\C}(\zeta) = 0 \}.
\end{equation}
There exists  a natural current of integration over the nodal
hypersurface in any ball bundle $B^*_{\epsilon} M$ with $\epsilon
< \epsilon_0$ , given by
\begin{equation} \langle [Z_{\phi_{\lambda}^{\C}}] , \phi \rangle =  \frac{i}{2 \pi} \int_{B^*_{\epsilon} M} \ddbar \log
|\phi_{\lambda}^{\C}|^2 \wedge \phi =
\int_{Z_{\phi_{\lambda}^{\C}} } \phi,\;\;\; \phi \in \dcal^{ (m-1,
m-1)} (B^*_{\epsilon} M). \end{equation} In the second equality we
used the Poincar\'e-Lelong formula. The notation $\dcal^{ (m-1,
m-1)} (B^*_{\epsilon} M)$ stands for smooth test $(m-1,
m-1)$-forms with support in $B^*_{\epsilon} M.$

The nodal hypersurface $Z_{\phi_{\lambda}^{\C}}$ also carries a
natural volume form $|Z_{\phi_{\lambda}^{\C}}|$ as a complex
hypersurface in a \kahler manifold. By Wirtinger's formula, it
equals the restriction of $\frac{\omega_g^{m-1}}{(m - 1)!}$ to
$Z_{\phi_{\lambda}^{\C}}$. Hence, one can regard
$Z_{\phi_{\lambda}^{\C}}$ as defining  the  measure
\begin{equation} \langle |Z_{\phi_{\lambda}^{\C}}| , \phi \rangle
= \int_{Z_{\phi_{\lambda}^{\C}} } \phi \frac{\omega_g^{m-1}}{(m -
1)!},\;\;\; \phi \in C(B^*_{\epsilon} M).
\end{equation}
 We prefer to state results in terms of the
current $[Z_{\phi_{\lambda}^{\C}}]$ since it carries more
information.

We will say that a sequence $\{\phi_{j_k}\}$ of $L^2$-normalized
eigenfunctions  is {\it quantum ergodic} if
\begin{equation} \label{QEDEF} \langle A \phi_{j_k}, \phi_{j_k} \rangle \to
\frac{1}{\mu(S^*M)} \int_{S^*M} \sigma_A d\mu,\;\;\; \forall A \in
\Psi^0(M). \end{equation} Here, $\Psi^s(M)$ denotes the space of
pseudodifferential operators of order $s$, and $d \mu$ denotes
Liouville measure on the unit cosphere bundle $S^*M$ of $(M, g)$.
More generally, we denote by $d \mu_{r}$ the (surface) Liouville
measure on $\partial B^*_{r} M$, defined by
\begin{equation} \label{LIOUVILLE} d \mu_r = \frac{\omega^m}{d |\xi|_g} \;\; \mbox{on}\;\; \partial B^*_r
M.
\end{equation}

Our main result is:

\begin{theo}\label{ZERO}  Let $(M, g)$ be  real analytic, and let $\{\phi_{j_k}\}$ denote a quantum ergodic sequence
of eigenfunctions of its Laplacian $\Delta$.  Let
$(B^*_{\epsilon_0} M, J)$ be the maximal Grauert tube around $M$
with complex structure $J_g$ adapted to $g$. Let $\epsilon <
\epsilon_0$. Then:
$$\frac{1}{\lambda_j} [Z_{\phi_{j_k}^{\C}}] \to  \frac{i}{ \pi} \overline{\partial} {\partial} |\xi|_g
 = \frac{1}{2 \pi |\xi|_g}
 \omega_g + \frac{ d |\xi|_g^2 \wedge \alpha}{4 \pi  |\xi|_g^3},\;\;
 \mbox{weakly in }\;\;  \dcal^{' (1,1)} (B^*_{\epsilon} M). $$
\end{theo}

In other words, for any continuous test form $\psi \in \dcal^{'
(m-1, m-1)}(B^*_{\epsilon} M)$, we have
$$\frac{1}{\lambda_j} \int_{Z_{\phi_{j_k}^{\C}}} \psi \to
 \frac{i}{ \pi} \int_{B^*_{\epsilon} M} \psi \wedge \overline{\partial} {\partial}
|\xi|_g. $$ The  limit
 $(1,1)$ form $\frac{i}{ \pi} \overline{\partial} {\partial} |\xi|_g$
 first arose  in \cite{LS1, GS1} for a different reason which
 is reviewed in
 \S \ref{EXAMPLES}.
As a corollary we obtain a similar result on the integrals of
scalar functions against  the measures $|Z_{\phi_{j_k}^{\C}}|$:
for any $\phi \in C(B^*_{\epsilon} M)$,
$$\frac{1}{\lambda_j} \int_{Z_{\phi_{j_k}^{\C}}} \phi \frac{\omega_g^{m-1}}{(m -
1)!}  \to
 \frac{i}{ \pi} \int_{B^*_{\epsilon} M} \phi  \overline{\partial} {\partial}
|\xi|_g \wedge \frac{\omega_g^{m-1}}{(m - 1)!} . $$

As is well-known, ergodic sequences of density one in the spectrum
arise when the geodesic flow is ergodic, and an entire orthonormal
basis is ergodic when the Laplacian is quantum uniquely ergodic
\cite{CV, Shn, Z, Z2}. Thus, we obtain the titled result:

\begin{cor}\label{ZEROCOR}  Let $(M, g)$ be a real analytic with ergodic  geodesic
flow.  Let $\{\phi_{j_k}\}$ denote a full density ergodic
sequence. Then for all $\epsilon < \epsilon_0$,
$$\frac{1}{\lambda_{j_k}} [Z_{\phi_{j_k}^{\C}}] \to  \frac{i}{ \pi} \overline{\partial} {\partial} |\xi|_g,\;\;
 \mbox{weakly in}\;\; \dcal^{' (1,1)} (B^*_{\epsilon} M). $$
\end{cor}

By the unique quantum ergodicity result of E. Lindenstrauss
\cite{Lin}, the zero currents of the full sequence of Hecke
eigenfunctions on arithmetic hyperbolic surfaces satisfy the
 limit formula in Corollary \ref{ZEROCOR}. The Rudnick-Sarnak conjecture \cite{RS} that
negatively curved compact manifolds are quantum uniquely ergodic
would imply that the zero currents of the full sequence of
eigenfunctions should satisfy the limit formula on such spaces.

Theorem \ref{ZERO} also has implications for a general Riemannian
manifold $(M, g)$. The orthonormal basis of a general Riemannian
manifold is not quantum ergodic and moreover the complex zeros do
not generally tend to the limit $\frac{i}{ \pi}
\overline{\partial} {\partial} |\xi|_g$ (a flat torus provides a
simple example where zeros concentrate on complex hypersurface).
However, in a precise sense, a {\it random orthonormal basis} of
$L^2(M)$ (adapted to $\Delta_g$) has the quantum ergodic property,
and hence the complex zeros of the basis functions will satisfy
the limit formula of Theorem \ref{ZERO}.

To state the result, we  recall the definition and results on
these random orthonormal bases from \cite{Z3}. We partition the
spectrum of $\sqrt{\Delta_g}$ into the intervals $I_k = [k, k +
1]$ and denote by $\Pi_k = E(k+1) - E(k)$ the spectral projections
for $\sqrt{\Delta_g}$ corresponding to the interval $I_k$. We
denote by $N(k)$ the number of eigenvalues
  in $I_k$ and  put $\hcal_k = \mbox{ran} \Pi_k$ (the range of
 $\Pi_k$). $\hcal_k$ consists of
 linear combinations $\sum_{j: \lambda_j \in I_k} c_j
\phi_j$  of the eigenfunctions of $\sqrt{\Delta_g}$ with
eigenvalues in $I_k$.  We define
 a  {\it random} orthonormal basis $\{U_k \phi_j\}$ of $\hcal_k$ by
changing the basis of $\sqrt{\Delta}$-eigenfunctions $\{\phi_j\}$
of $\Delta$ in $\hcal_k$ by a random element $U_k$ of the unitary
group $U(\hcal_k)$ (equipped with its normalized Haar measure $d
\nu_k$ of the finite dimensional Hilbert space $\hcal_k$.

 We
then define a {\it random orthonormal basis} of $L^2(M)$ (adapted
to $\Delta_g$) by taking the product over all the spectral
intervals in our partition. That is,  we define the infinite
dimensional unitary group
$$U(\infty) =\Pi^{\infty}_{k=1}  U(\hcal_k)$$
 of  sequences
 $(U_1, U_2,\dots)$, with $U_k \in U(\hcal_k)$, and
 equip $U(\infty)$ with the
 product measure
$$d\nu_{\infty}= \Pi^{\infty}_{k=1} d\nu_k.$$
A  {\it random} orthonormal basis $\Psi = \{ (U_k \phi_j) \}$ of
$L^2(M)$ is thus an orthonormal basis obtained by applying a
random element $U \in U(\infty)$ to the orthonormal basis $\Phi =
\{\phi_j\}$ of eigenfunctions of $\sqrt{\Delta}$.

In  \cite{Z3}, it is proved  that  random orthonormal bases
satisfy the following variance asymptotics:
$$\begin{array}{l}   \E \left(
 \sum_{j: \lambda_j \in I_k}
|(AU\phi_j, U\phi_j)-\omega(A)|^2 \right) \sim  ( \omega(A^*A) -
\omega(A)^2 ).
\end{array}$$
To be precise, in \cite{Z3} it is assumed that the widths of the
intervals $I_k$ increase to infinity, but    the more recent
strong Szeg\"o asymptotics of \cite{GO, LRS} allow one to prove
the same result for intervals of bounded width such as the $I_k$
above. By the strong law of large numbers (see \cite{Z3}) it
follows that with probability one, a random orthonormal basis of
$L^2(M)$ is quantum ergodic. We thus have:

\begin{cor}\label{ZERORAN}  Let $(M, g)$ be any real analytic compact
Riemannian manifold. Then with probability one, a random
orthonormal basis $\{ \psi_j = U \phi_j\}$ of $L^2(M)$ as defined
above satisfies
$$\frac{1}{\lambda_{j_k}} [Z_{\psi_{j_k}^{\C}}] \to  \frac{i}{ \pi} \overline{\partial} {\partial} |\xi|_g,\;\;
 \mbox{weakly in}\;\; \dcal^{' (1,1)} (B^*_{\epsilon} M), $$
 for a full density  subsequence $\{\psi_{j_k}\}$.
\end{cor}
This  gives a kind of almost sure  improvement of the
complexification of Theorem 14.3 of Jerison-Lebeau \cite{JL} from
an inequality to an asymptotic formula.

\subsection{\label{DISCUSSION} Discussion and outline of the proof}

In summary,  the complex zeros of a quantum ergodic sequence
become equidistributed with respect to the $(1,1)$ form $
\frac{i}{\pi} \overline{\partial} {
\partial} |\xi|_g$.
As mentioned above, the  \kahler form  on $B^*_{\epsilon_0} M$
associated to $g$ is the (1,1)-form $\omega_g = \frac{1}{i}
\overline{\partial}
\partial |\xi|_g^2$.
  We observe that  $ \overline{\partial}
{\partial} |\xi|_g$  is singular relative to the \kahler form
along the zero section, i.e. the totally real submanifold $M$ (the
geometry will be reviewed in \S \ref{AC}). This singular
concentration could be attributed to the fact that the Laplacian
is time reversal invariant (i.e. invariant under complex
conjugation), so that the eigenfunctions are usually
 real-valued on the real
submanifold $M$. Hence their complex zero sets are invariant under
the (time reversal) involution $\sigma: (x, \xi) \to (x, -\xi)$
(the classical limit of complex conjugation).

As a simple example of Theorem \ref{ZERO}, consider the circle
$S^1$. The geodesic flow is  ergodic modulo the symmetry $(x, \xi)
\to (x, - \xi)$. The real eigenfunctions $\sin 2 \pi k x, \cos 2
\pi k x$ are therefore quantum ergodic. They complexify to the
cylinder as $\sin 2 \pi k z, \cos 2 \pi k z$. The complex zero set
of these holomorphic functions lies entirely on the set $\Im z =
0$ and become uniformly distributed with respect to $2 d \theta$
as $k \to \infty$. It will be checked in \S \ref{PROOF} that the
coefficient agrees with the result of Theorem \ref{ERGO}. Note
however that  the complex eigenfunctions $e^{\pm 2 \pi i k x}$ are
not quantum ergodic and have no complex zeros (see \S \ref{PROOF}
for further discussion).

We now outline the  main steps of the proof. As mentioned above,
eigenfunctions $\phi_{\lambda}$ of Laplcians of  real analytic
Riemannian manifolds admit holomorphic extensions
$\phi_{\lambda}^{\C}$ to a maximal Grauert tube in the
complexification $M_{\C}$ of $M$, which we will identify with a
maximal ball bundle  $B^*_{\epsilon_0} M$ on which the adapted
complex structure is defined. The square $|\xi|_g^2$  of  metric
norm function $|\xi|_g = \sqrt{\sum_{i, j = 1}^m g^{ij}(x) \xi_i
\xi_j}$ is a smooth, strictly plurisubharmonic exhaustion function
on $B^*_{\epsilon_0} M$.  For $0 < \epsilon \leq \epsilon_0$ the
sphere bundles $S^*_{\epsilon} M =
\partial B^*_{\epsilon} M$ are strictly pseudoconvex CR manifolds.
We denote by $\ocal(B_{\epsilon}^*(M))$ the class of holomorphic
functions on this domain, and by $\ocal(\partial
B_{\epsilon}^*(M))$ the space of boundary values of holomorphic
functions, i.e. the CR holomorphic functions.  For each $0<
\epsilon < \epsilon_0$, the  restriction
$\phi_{\lambda}^{\C}|_{\partial B_{\epsilon}^*(M)}$ thus lies in
the Hardy space $\ocal^0(\partial B_{\epsilon}^*(M))$ of square
integrable CR functions.

A key object in the proof  is the sequence of functions
$U_{\lambda}(x, \xi) \in C^{\infty}( B^*_{\epsilon} M)$ defined by
\begin{equation} \label{DEFS} \left\{ \begin{array}{l}  U_{\lambda}(x, \xi) : =
\frac{\phi_{\lambda}^{\C}(x, \xi)}{\rho_{\lambda}(x, \xi)},\;\;\;
(x, \xi) \in
 B^*_{\epsilon} M, \;\;\; \mbox{where}\\ \\
\rho_{\lambda} (x, \xi) :=  ||\phi_{\lambda}^{\C} |_{\partial
B_{|\xi|_g} } ||_{L^2(\partial B^*_{|\xi|_g}M)}
\end{array} \right.
\end{equation}
Thus,  $\rho_{\lambda}(x, \xi)$ is the `moving'   $L^2$-norm of
$\phi_{\lambda}^{\C}$ as it is restricted to the one-parameter
family $\{\partial B^*_{\epsilon } M\}$ of strictly pseudo-convex
CR manifolds.  $U_{\lambda}$ is of course not  holomorphic, but
its restriction  to each sphere bundle is CR holomorphic there,
i.e.
\begin{equation} \label{LITTLEU} u_{\lambda}^{\epsilon}  =
U_{\lambda} |_{\partial B^*_{\epsilon} M}  \in \ocal^0 (\partial
B^*_{\epsilon}(M).
\end{equation}

Our first result gives an ergodicity property of holomorphic
continuations of ergodic eigenfunctions.

\begin{lem} \label{ERGOCOR} Assume that $\{\phi_{\lambda}\}$ is a quantum ergodic sequence of $\Delta$-eigenfunctions
on $M$ in the sense of (\ref{QEDEF}).  Then for each $0 < \epsilon
<\epsilon_0$,
$$|U_{\lambda}|^2  \to  \frac{1}{\mu_1(S^* M)} |\xi|_g^{-m + 1},\;\; \mbox{weakly in}\;\; L^1(
B^*_{\epsilon}M, \omega^m) . $$
\end{lem}

We note that $\omega^m = r^{m-1} dr d\omega dvol(x)$ in polar
coordinates, so the right side indeed lies in $L^1$. The actual
limit function is otherwise irrelevant. The next step is to use a
compactness argument from \cite{SZ} (see also \cite{NV}) to obtain
strong convergence of the normalized logarithms of the sequence
$\{|U_{\lambda}|^2\}$. The first statement of the following lemma
immediately implies the second.

\begin{lem} \label{ZEROWEAK} Assume that $|U_{\lambda}|^2  \to  \frac{1}{\mu_1(S^* M)}  |\xi|_g^{-m + 1
},\;\; \mbox{weakly in}\;\; L^1( B^*_{\epsilon}M, \omega^m). $
Then:

\begin{enumerate}

\item  $\frac{1}{\lambda_j}  \log |U_j|^2 \to
0$ strongly in $L^1(B^*_{\epsilon} M). $

\item $\frac{1}{\lambda_j} \ddbar \log |U_j|^2 \to 0,\;\;
\mbox{weakly in} \;\; \dcal{'(1,1)}(B^*_{\epsilon} M). $

\end{enumerate} \end{lem}

Separating out the numerator and denominator of $|U_j|^2$, we
obtain that
 \begin{equation} \label{SEPARATE} \frac{1}{\lambda_j} \ddbar  \log |\phi_{\lambda}^{\C}|^2 - \frac{2}{\lambda_j} \ddbar
 \log \rho_{\lambda_j} \to 0,\;\;\; (\lambda_j \to \infty). \end{equation}  The next lemma shows that the second term
has a weak limit:

\begin{lem} \label{NORM}  For $0 < \epsilon < \epsilon_0$,
$$\frac{1}{\lambda} \log \rho_{\lambda}(x, \xi)  \to \; |\xi|_{g_x},\;\;\;\mbox{ in }\;\;L^1(B^*_{\epsilon} M)\;\; \mbox{as}\;\; \lambda \to
\infty.
$$
Hence, $$\frac{1}{\lambda_j} \ddbar
 \log \rho_{\lambda_j} \to \ddbar |\xi|_{g_x},\;\;\; (\lambda_j \to \infty) \;\; \mbox{weakly in }\;\;  \dcal'(B^*_{\epsilon}
M).$$

\end{lem}

It follows that the left side of (\ref{SEPARATE}) has the same
limit, and that will complete the proof of
 Theorem \ref{ZERO}.

The proofs of the lemmas are based on the properties of the
analytic continuation of the wave kernel  as  a complex Fourier
integral associated to the complexified exponential map
 \cite{Bou, GS2, GLS}.   Since
 $\phi_{\lambda}^{\C} |_{\partial B^*_{ \sqrt{\rho}} M}  =
e^{- \lambda \sqrt{\rho}} E(i \sqrt{\rho}) \phi_{\lambda}$, the
analytic continuation of $\phi_{\lambda}$ is obtained by applying
a complex Fourier integral operator of known order and symbol.
This allows us to connect growth and distribution of zeros in
tubes  to the dynamics of the geodesic flow. Although this paper
treats only the ergodic case, it would be of interest to
investigate  complex zeros of $\Delta$-eigenfunctions under other
dynamical hypotheses such as complete integrability. It would also
be interesting to investigate analogues  for boundary value
problems.

Let us compare the results of this paper to earlier results of B.
Shiffman-S. Zelditch \cite{SZ} and of S.Nonnemacher- A. Voros
\cite{NV} (see also \cite{R}) on complex zeros of eigenfunctions
of ergodic quantum maps in \kahler phase spaces. These articles
were concerned with the complex zeros of eigenfunctions of ergodic
quantum maps acting on spaces $H^0(M, L^N)$ of sections of powers
of a holomorphic line bundle over a
 a \kahler manifold $(M, \omega)$. The role of the Grauert tube was played by the disc bundle $D^*
\subset L^*$ of the dual line bundle, and the role of the CR
manifolds $S^*_{\epsilon} M$ was played by the circle bundle $X =
\partial D^*$, which is a strictly pseudoconvex CR submanifold of
$L^*$. The norm function $|\xi|_g$ of $(M, g)$ is thus  analogous
to the hermitian metric  of $(L, h)$, but the analogy is not very
close.  Indeed, the role of the geodesic role in the $(M, g)$
setting was split in the \kahler setting between two dynamical
systems: the $S^1$ action on the circle bundle $X \to M$ (which is
the direct analogue of the geodesic flow in the Riemannian setting
but is of course not ergodic); and an auxiliary ergodic symplectic
transformation $\chi$ of the \kahler manifold $M$. Ergodicity of
joint eigenfunctions $\{s_{N, j}\}$ of the $S^1$ action and of the
quantum map  associated to $\chi$  gave $||s_{N,j}(z)||_{h^N}^2
\to 1$ in the weak sense (as in Theorem \ref{ERGOCOR}), where the
norm is the pointwise norm relative to a hermitian metric on $L$
with curvature equal to the \kahler form $\omega$. Weak
convergence of the zero currents as in Lemma \ref{ZEROWEAK} showed
that $\frac{1}{N} Z_{s_{N,j}} \to \omega,$ proving
equidistribution of zeros relative to the \kahler form. In
comparison, Theorem \ref{ZERO} shows that complex zeros of real
ergodic eigenfunctions are not equidistributed relative to the
analogous \kahler form $\omega_g$ in the Riemannian setting, but
rather to the relatively singular form stated in the theorem. The
difference can be traced to Lemma \ref{NORM}, which indicates that
the moving $L^2$-norm $\rho_{\lambda}$ ends up playing the key
role of the hermitian metric.

\subsection{\label{CONJ} A conjecture on ergodic real nodal hypersurfaces}

We close the introduction by stating a conjecture on the real
nodal hypersurfaces $Z_{\phi_j} = \{x \in M: \phi_j(x) = 0\}$  on
Riemannian manifolds with ergodic geodesic flow. We define the
distribution of real zeros of $\phi_{\lambda}$ by integration
\begin{equation} \label{REALINT} \langle [Z_{\phi_j}], f \rangle =
\int_{Z_{\phi_j}} f(x) d {\mathcal H}^{m-1}, \end{equation} with
respect to    $(m-1)$-dimensional Haussdorf measure $d{\mathcal
H}^{m-1}$  on the nodal hypersurface induced by the Riemannian
metric of $(M, g)$.

\begin{conj} \label{REALZ} Let $(M, g)$ be a real analytic Riemannian manifold
with ergodic geodesic flow, and let $\{\phi_j\}$ be the density
one sequence of ergodic eigenfunctions. Then,
$$\langle [Z_{\phi_j}], f \rangle \sim  \{ \int_M f
dVol_g\} \lambda. $$
\end{conj}

More generally, we  conjecture the same limit result for any
quantum ergodic sequence of eigenfunctions. In the case of random
spherical harmonics, the limit formula was proved in the PhD
thesis of J. Neuheisel \cite{Ne}. We believe that (in a
straightforward way) it can be extended to random orthonormal
bases on any compact Riemannian manifold as defined above, which
would give the asymptotic strengthening of \cite{JL} in the real
domain.

 There is of course a very wide gap between the equidistribution  Conjecture
\ref{REALZ} and the best known result on volumes (\ref{DF}). Due
to the singular concentration along the real zero set, it seems
possible that Corollary \ref{ERGOCOR} could have implications for
the distribution of real zeros in the ergodic case.
\medskip

\noindent{\bf Acknowledgements} This work developed out of joint
work with B. Shiffman on line bundles. A  preliminary version  was
presented at the Newton Institute workshop on quantum chaos in
June, 2004 and the final version was complete at the IHP program
Time at Work in June 2005. We thank these institutions for their
support.  We would also like to thank Z. Rudnick, M. Sodin for
suggesting improvements in the exposition and D. Jerison for
encouragement to include Corollary \ref{ZERORAN}.

\section{\label{AC} Analytic continuation  to a Grauert tube}

In this section, we recall the relevant known results on analytic
continuation of the wave kernel and $\Delta$- eigenfunctions of a
real analytic $(M, g)$ to a Grauert tube.

A real analytic manifold $M$ always possesses a complexification
$M_{\C}$, i.e. a complex manifold  of which $M$ is a totally real
submanifold (Bruhat-Whitney \cite{BW}). The germ of $M_{\C}$ along
$M$ is unique. In \cite{G}, Grauert constructed plurisubharmonic
exhaustion functions $\rho$ on $M_{\C}$, which define the Grauert
tubes $M_{\epsilon} = \{\rho < \epsilon\}$ relative to $\rho.$

In \cite{GS1, GS2, LS1, LS2, GLS}, the complex geometry of Grauert
tubes was brought into contact with the symplectic geometry of the
cotangent bundle $T^*M$ and with the Riemannian geometry of real
analytic metrics $g$. The key results (for the purposes of this
article) are the following: First, a metric $g$ determines a
canonical plurisubharmonic function $\rho_g$ on $M_{\C}$. It is
defined on a maximal Grauert tube whose radius $\epsilon_0$  is
often called the `radius of the Grauert tube'.  There is a
symplectic diffeomorphism
$$\psi: M_{\epsilon} \to B^*_{\epsilon} M$$ of the Grauert tubes
with respect to $\rho_g$ to the ball bundles $B^*_{\epsilon} M
\subset T^*M$ with respect to $g$, which identifies $\rho_g$ with
$|\xi|_g^2$. As pointed out in \cite{LS1} (see also \cite{GLS}),
one may take $\psi^{-1}$ to be   the complexified exponential map
$$(x, \xi) \in B^*M \to \exp_x \sqrt{-1} \xi \in M_{\epsilon}.$$
The map $\psi$ endows $B^*_{\epsilon} M$ with a complex structure
$J_g$  {\it adapted} to $g$.

\begin{prop} \label{EXPMAP} \cite{LS1,GLS} The adapted
complex structure is uniquely characterized by the property that
the complexified exponential map,
$$(x, \xi) \in B_{\epsilon}^*M \to \exp_x \sqrt{-1} \xi \in M_{\epsilon}$$
is a biholmorphism for $0 < \epsilon < \epsilon_0$. \end{prop}

The domains $B^*_{\epsilon} M$ are strictly pseudoconvex for
$\epsilon < \epsilon_0$,  hence their boundaries $S^*_{\epsilon} M
= \partial B^*_{\epsilon} M$ are strictly pseudoconvex CR
manifolds. We use the notation $\partial B^*_{\epsilon} M$ to
emphasize the role of sphere bundles as boundaries of domains. The
restrictions
$$\psi_{\epsilon} = \exp (i \epsilon)^{-1}: \partial M_{\epsilon} \to \partial
B^*_{\epsilon}$$ of $\psi$ are then CR holomorphic
diffeomorphisms.

As mentioned above, the metric norm function $|\xi|_g$ pulls back
under $\psi_{\epsilon}$ to the function
 $  \sqrt{\rho}$ on $M_{\C}$, which is
known as the  Monge-Amp\'ere function \cite{GS1}. It  equals
$r_{\C}(z, \bar{z})$ where $r_{\C}$ is the holomorphic extension
of the distance function.  In the cotangent picture, the metric
norm function $|\xi|_g$ is smooth on $B^*_{\epsilon_0}M \backslash
M$ and solves the homogeneous complex Monge-Amp\`ere equation
$(\ddbar |\xi|_g)^m = 0$ there. In fact, the form $\ddbar |\xi|_g$
has rank $m - 1$ on $B^*_{\epsilon_0}M \backslash M$, and its
kernel is a smooth rank $1$ sub-bundle of $T (B^*_{\epsilon_0}M
\backslash M)$. The leaves of the associated (`Monge-Amp\`ere' or
Riemann) foliation are the complex curves $t + i \tau \to \tau
\dot{\gamma}(t)$, where $\gamma$ is a geodesic, where $\tau > 0$
and where  $\tau \dot{\gamma}(t)$ denotes multiplication of the
tangent vector to $\gamma$ by $\tau$. We refer to \cite{LS1} for
further discussion.

\subsection{\label{EXAMPLES} Model examples}

To better understand complexifications and Monge-Amp\'ere
functions, and in particular the limit form in Theorem \ref{ZERO},
we go over  several (well-known) model examples. We note that
these examples do not have ergodic geodesic flow, so they do not
exemplify Theorem \ref{ZERO}, but only the objects involved in it.
\medskip

\noindent{\bf (i) } Complex tori:
\medskip

The complexification of the  torus $M = \R^m/\Z^m$ is $M_{\C} =
\C^m/\Z^m$. The adapted complex structure to the flat metric on
$M$  is the standard (unique)  complex structure on $\C^m$. The
complexified exponential map is $\exp_{\C x} (i \xi) = z: =  x + i
\xi $, while  the distance function $r(x, y) = |x - y|$ extends to
$r_{\C}(z, w) = \sqrt{(z - w)^2}. $ Then $\sqrt{\rho}(z, \bar{z})
= \sqrt{(z - \bar{z})^2} = \pm 2 i |\Im z| = \pm 2 i |\xi|. $
Thus, the limit form is $\frac{i}{\pi} \ddbar |\Im z|.$

\noindent{\bf  (ii)  $\Ss^n$} \cite{PW, GS1}  The unit sphere
$x_1^2 + \cdots + x_{n+1}^2 = 1$ in $\R^{n+1}$ is complexified as
the complex quadric
$$S^2_{\C} = \{(z_1, \dots, z_n) \in \C^{n + 1}: z_1^2 + \cdots +
z_{n+1}^2 = 1\}. $$ If we write $z_j = x_j + i \xi_j$, the
equations become $|x|^2 - |\xi|^1 = 1, \langle x, \xi \rangle =
0$. The geodesic flow $G^t(x, \xi) = (\cos t x + \sin t \xi, -
\sin t x + \cos t \xi)$  induces the exponential map $exp_x  \xi =
(\cos |\xi| ) x + (\sin |\xi|) \xi, $ which complexifies to
$$\exp_{\C, x} \sqrt{-1} \xi = (\cosh |\xi|) x +  \sqrt{-1} (\sinh |\xi|)\frac{ \xi}{|\xi|}. $$

The distance function of $\Ss^n$ of constant curvature $1$ is
given by:
$$r(x, y) = 2 \sin^{-1} \frac{|x - y|}{2} = 2 \sin^{-1}
(\frac{1}{2} \sqrt{(x - y)^2}), $$ whose analytic continuation to
$\Ss^n_{\C} \times \Ss^n_{\C}$ is the doubly-branched holomorphic
function:
$$r_{\C}(z, w) = 2 \sin^{-1} \frac{1}{2} \sqrt{(z - w)^2}. $$
One branch gives the pluri-subharmonic function
$$\sqrt{\rho}(z) = r_{\C}(z, \bar{z}) = 2 \sin^{-1} i |\Im z|) = 2 i \sinh^{-1}  |\Im z|) =  i \cosh^{-1} |z|^2, \;\;(z \in \Ss^n_{\C}). $$
 Since
$$\begin{array}{lll} \exp_{\C}^* \sqrt{\rho } (x, \xi) & = & \cosh^{-1} |
(\cosh |\xi|) x +  i (\sinh |\xi|)\frac{ \xi}{|\xi|}|^2
\\ & = & \cosh^{-1} \{(\cosh |\xi|)^2 - (\sinh |\xi|)^2\} \\
& = & \cosh^{-1}  \cosh 2 |\xi| = 2 |\xi|, \end{array}$$ the limit
form is $\frac{i}{2 \pi} \ddbar  \cosh^{-1} |z|^2$ on
$\Ss^n_{\C}$.
\medskip

\noindent{\bf (iii) (See e.g. \cite{KM}).  $\Hh^n$} The
hyperboloid model of hyperbolic space is the hypersurface in
$\R^{n+1}$ defined by
$$\Hh^n = \{ x_1^2 + \cdots
x_n^2 - x_{n+1}^2 = -1, \;\; x_n > 0\}. $$ Then,
$$H^n_{\C} = \{(z_1, \dots, z_{n+1}) \in \C^{n+1}:  z_1^2 + \cdots
z_n^2 - z_{n+1}^2 = -1\}. $$ In real coordinates $z_j = x_j + i
\xi_j$, this is:
$$\langle x, x \rangle_L - \langle \xi, \xi \rangle_L = -1,\;\;
\langle x, \xi \rangle_L = 0$$ where $\langle, \rangle_L$ is the
Lorentz inner product of signature $(n, 1)$. The complexified
exponential map is given by $$\exp_{\C x} (\sqrt{-1} \xi) = \cos
(\frac{||\xi||_L}{\sqrt{2}}) x + \sqrt{-1} (\frac{ \sin
\frac{||\xi||_L}{\sqrt{2}}}{||\xi||_L}) \xi.
$$
Let
$$\begin{array}{l} M_{\epsilon} =  \{z \in \C^{n+1}: z_1^2 + \cdots + z_n^2 - z_{n+1}^2 = -1,\;\; |z_1|^2 + \cdots + |z_n|^2 - |z_n|^2 < \epsilon \}.\end{array}$$
We note that $M_{\epsilon}$ has two components according to the
sign of $\Re z_{n+1}.$ The Monge-Ampere function is:
$$\sqrt{\rho}(z) = \cos^{-1} (||x||_L^2 + ||\xi||_L^2 - \pi)/\sqrt{2}. $$
 The radius of maximal Grauert
tube is $\epsilon = 1$ or $r = \pi/\sqrt{2}.$ Hence the limit form
is $\frac{i}{\pi} \ddbar \cos^{-1} (||x||_L^2 + ||\xi||_L^2 -
\pi)/\sqrt{2})$ on $M_1$.

\subsection{Analytic Continuation of the wave kernel}

By the wave kernel of $(M, g)$ we mean the kernel $$E(t, x, y) =
\sum_{j = 0}^{\infty} e^{i t \lambda_j} \phi_j(x) \phi_j(y)$$ of
$e^{i t \sqrt{\Delta}}.$
 As discussed in
\cite{Bou, GS2, GLS}, the wave kernel at imaginary times admits a
holomorphic extension to $M_{\epsilon} \times M$ as
\begin{equation} \label{EI} E(i \epsilon, \zeta, y) = \sum_{j = 0}^{\infty} e^{-
\epsilon \lambda_j} \phi_{C j} (\zeta) \phi_j(y),\;\;\; (\zeta, y)
\in M_{\epsilon}  \times M.  \end{equation}

In the simplest (albeit non-compact) case of $\R^n$, the wave
kernel $E(t, x, y) = \int_{\R^n} e^{i t |\xi|} e^{i \langle \xi, x
- y \rangle} d\xi$ analytically continues  to $t + i \tau, \zeta =
x + i p \in \C_+ \times \C^n$ as the integral
$$E(t + i \tau , x + i p , y) = \int_{\R^n} e^{i (t + i \tau)  |\xi|} e^{i \langle \xi, x + i p - y
\rangle} d\xi, $$ which converges absolutely for $|p| < \tau.$ At
positive imaginary times and for $(x, y) \in \R^n \times \R^n$,
$E(i \tau, x, y)$ is the Poisson kernel of the upper half space
$\R^n_x \times \R_{\tau}^+$,
$$K(\tau, x, y) = \tau^{-n} \left(1 + (\frac{x -
y}{\tau})^2 \right)^{-\frac{n+1}{2}}  = \tau \left(\tau^2 + (x -
y)^2)\right)^{-\frac{n + 1}{2}},$$ which visibly has a holomorphic
continuation to $\zeta = x + i p$ in the $x$ variable for $|p| <
\tau$.

 On
a general analytic Riemannian manifold, one has a similar integral
formula for the wave kernel of the form \begin{equation}
\label{PARAONE} E(t, x, y) = \int_{T^*_y M} e^{i t |\xi|_{g_y} }
e^{i \langle \xi, \exp_y^{-1} (x) \rangle} A(t, x, y, \xi) d\xi
\end{equation} where $|\xi|_{g_x} $ is the metric norm function at
$x$, and where $A(t, x, y, \xi)$ is a polyhomogeneous amplitude of
order $0$. The analytic continuation of the wave group at
imaginary times is the Poisson operator $e^{- \tau
\sqrt{\Delta}}$, a Fourier integral operator with complex phase;
for background on the microlocal analysis of such kernels, we
refer to \cite{T} (Chapter XI). The analytic continuation,
\begin{equation} \label{CXPARAONE} E(i \tau,
\zeta, y) = \int_{T^*_y} e^{- \tau  |\xi|_{g_y} } e^{i \langle
\xi, \exp_y^{-1} (\zeta) \rangle} A(t, \zeta, y, \xi) d\xi
\;\;\;(\zeta = x + i p).
\end{equation} of the Poisson kernel to the Grauert tube $|\zeta|
< \tau$ defines a complex Fourier integral operator from $L^2(M)$
with values in holomorphic functions in $M_{\tau}$ where $\tau <
\epsilon_0$. Its  canonical relation is the complexification of
the canonical relation of the real wave group, i.e. `graph' of the
complexified geodesic flow at imaginary times
\begin{equation} \label{GAMMAIT} \Gamma_{i \tau} = \{(x, \xi, z, \zeta) \in T^*M \backslash 0
\times T^*M_{\epsilon} \backslash 0 : G^{i \tau}(x, \xi) = (z,
\zeta). \} \end{equation} We may (and will) also express the
adjoint operator in the modified form
\begin{equation} \label{CXPARAONEADJ} E^*(i \tau,
x, \zeta) = \int_{\R^n} e^{- \tau  |\xi|_{g_x} } e^{i \langle \xi,
\exp_x^{-1} (\zeta) \rangle} A^*(t, \zeta, x, \xi) d\xi
\;\;\;(\zeta = x + i p).
\end{equation}
Here, we use the phase $- \tau |\xi|_{g_x} + i \langle \xi,
\exp_x^{-1}(\bar{\zeta}) \rangle$ instead of $- \tau |\xi|_{g_x} -
i \langle \xi, \exp_{\bar{\zeta}}^{-1}(x) \rangle.$

We can transport the complexified wave kernel to $B^*_{\epsilon} M
\times M$ using the complexified exponential map (or $\psi$):
\begin{equation} \label{ETILDE}  \tilde{E}(i \epsilon, (x, \xi) ,
y) : = E(i \epsilon, \exp^{\C}_x \sqrt{-1} \xi, y). \end{equation}
It is again a complex Fourier integral operator whose canonical
relation $\tilde{\Gamma}_{i \tau}$ can be obtained from
(\ref{GAMMAIT}) by composing with $\psi$.

\subsection{Szego projector}

We  denote by $\ocal^{s + \frac{n-1}{4}}(\partial B^*_{\epsilon}
M)$ the subspace of the Sobolev space $W^{s +
\frac{n-1}{4}}(\partial B^*_{\epsilon} M)$ consisting of CR
holomorphic functions, i.e.
$${\mathcal O}^{s + \frac{m-1}{4}}(\partial B^*_{\epsilon} M) =
W^{s + \frac{m-1}{4}}(\partial B^*_{\epsilon} M) \cap \ocal
(\partial B^*_{\epsilon} M). $$ The inner product on $\ocal^0
(\partial B^*_{\epsilon} M)$ is with respect to the Liouville
measure $d\mu_{\epsilon}.$ There are similar spaces for $\partial
M_{\epsilon}$ and composition with $\psi_{\epsilon}$ defines an
isomorphism
$$\psi^*_{\epsilon}: {\mathcal O}^{s + \frac{m-1}{4}}(\partial B^*_{\epsilon}
M) \to {\mathcal O}^{s + \frac{m-1}{4}}(\partial M_{\epsilon}).
$$

We further  denote by $$ \tilde{\Pi}_{\epsilon} : L^2(\partial
B_{\epsilon}^* M) \to \ocal^0(
\partial B^*_{\epsilon} M)$$ the \szego projector for the tube $ B^*_{\epsilon}
M$, i.e. the orthogonal projection onto boundary values of
holomorphic functions in the tube. It is well-known  (cf.
\cite{BS, MS, GS2})  that $\tilde{\Pi}_{\epsilon}$ is a complex
Fourier integral operator, whose real canonical relation is the
graph $\Delta_{\Sigma}$  of the identity map on the symplectic
cone $\Sigma_{\epsilon} \subset T^* (\partial B^*_{\epsilon} M)$
spanned by the contact form $\alpha = \xi \cdot dx$, i.e.
$$\Sigma_{\epsilon} = \{(x, \xi; r \alpha_{\epsilon}),\;\;\;(x, \xi) \in \partial B^*_{\epsilon} M),\; r > 0 \}
 \subset T^* (\partial B^*_{\epsilon} M).\;\; $$
Alternatively it is a Toeplitz operator in the sense of Boutet de
Monvel-Guillemin \cite{BG}. The analogous \szego projector
$\Pi_{\epsilon}: L^2(\partial M_{\epsilon}) \to \ocal^0(
\partial M_{\epsilon})$ is conjugate to
$\tilde{\Pi}_{\epsilon}$ in the sense that $\tilde{\Pi}_{\epsilon}
= (\psi^*_{\epsilon})^{-1} \Pi_{\epsilon} \psi_{\epsilon}^*.$

We now consider the restrictions  $\tilde{\Pi}_{\epsilon} \circ
\tilde{E}(i \epsilon)$  from $L^2(M)$ to $\ocal (\partial
B^*_{\epsilon} M)$. Since $\Sigma_{\epsilon}$ is an $\R_+$-bundle
over $ \partial B^*_{\epsilon} M$, we can define the symplectic
equivalence of cones:
\begin{equation} \iota_{\epsilon} : T^*M \to
\Sigma_{\epsilon},\;\; \iota (x, \xi) = (x, \epsilon \xi, |\xi|
\alpha_{(x, \epsilon \xi)} ). \end{equation} The following result
is the transport under $\psi_{\epsilon}$ to $\partial
B^*_{\epsilon} M$ of results due to Boutet de Monvel (see also
\cite{GS2}).

\begin{theo}\label{BOUFIO}  \cite{Bou, GS2} \label{BDM} $\tilde{\Pi}_{\epsilon} \circ \tilde{E}(i \epsilon): L^2(M)
\to \ocal(\partial B^*_{\epsilon} M)$ is a  complex Fourier
integral operator of order $- \frac{m-1}{4}$  associated to the
canonical relation
$$\Gamma = \{(y, \eta, \iota_{\epsilon} (y, \eta) \} \subset T^*M \times \Sigma_{\epsilon}.$$
Moreover, for any $s$,
$$\tilde{\Pi}_{\epsilon} \circ \tilde{E} (i \epsilon): W^s(M) \to {\mathcal O}^{s +
\frac{m-1}{4}}(\partial B^*_{\epsilon} M)$$ is a continuous
isomorphism.
\end{theo}

\subsection{Analytic continuation of eigenfunctions}

We obtain the  holomorphic extension of the eigenfunctions
$\phi_{\lambda}$ by applying the complex Fourier integral operator
$E(i \tau)$:
\begin{equation} E(i \tau) \phi_{\lambda} = e^{- \tau \lambda}
\phi_{\lambda}^{\C}. \end{equation}  As usual, we can use the
complexified exponential map to transport the
$\phi_{\lambda}^{\C}$ to $B^*_{\epsilon_0}(M)$: \begin{equation}
\tilde{\phi}^{\C}_{\lambda} (x, \xi) = \phi_{\lambda}^{\C} (\exp_x
\sqrt{-1} \xi). \end{equation} By Theorem \ref{BOUFIO},we obtain:

\begin{cor} \cite{Bou, GLS} Each eigenfunction $\phi_{\lambda}$ has a
holomorphic extension to $ B^*_{\epsilon} M$ satisfying
$$\sup_{(x, \xi) \in B^*_{\epsilon}M} |\tilde{\phi}^{\C}_{\lambda}(x, \xi)| \leq
C_{\epsilon} \lambda^{m + 1} e^{\epsilon \lambda}. $$
\end{cor}

The fact that the holomorphic continuations of eigenfunctions can
be obtained by applying a complex Fourier integral operator is the
crucial link connecting the geodesic flow and the growth rate and
zeros of $\phi_{\lambda}^{\C}.$

\subsection{Examples}

We pause to consider some basic examples of holomorphic
continuations of eigenfunctions. The simplest  example is the flat
torus $\R^m/\Z^m,$  where the real eigenfunctions are $\cos
\langle k, x \rangle, \sin \langle k, x \rangle$ with $k \in 2 \pi
\Z^m.$ The complexified torus is $\C^m/\Z^m$ and the complexified
eigenfunctions are $\cos \langle k, \zeta \rangle, \sin \langle k,
\zeta \rangle$ with $\zeta  = x + i \xi.$

No such explicit examples exist in the ergodic case, but one can
see the uniform analytic continuation of eigenfunctions very
clearly  in the case  of compact hyberbolic quotients ${\bf
H}^m/\Gamma$. For simplicity we consider the two-dimensional case.
Eigenfunctions can be then represented by Helgason's generalized
Poisson integral formula \cite{H},
$$\phi_{\lambda}(z) = \int_B e^{(i \lambda + 1)\langle z, b
\rangle } dT_{\lambda}(b). $$ Here, $z \in D$ (the unit disc),
$B =
\partial D$, and $dT_{\lambda}
\in \dcal'(B)$ is the boundary value of $\phi_{\lambda}$, taken in
a weak sense along circles centered at the origin $0$. Also,
$\langle z, b \rangle$ is the (signed) hyperbolic distance of the
horocycle passing through $z$ and $b$ to $0$.
 To  analytically continue $\phi_{\lambda}$ it suffices  to
analytically continue $\langle z, b\rangle. $ Writing the latter
as  $\langle \zeta, b \rangle,  $ we have:
$$\phi_{\lambda}^{\C} (\zeta) = \int_B e^{(i \lambda + 1)\langle \zeta, b
\rangle } dT_{\lambda}(b). $$ Using this representation, one could
verify directly the estimates on growth of  complexified
eigenfunctions.

\section{The operators $ \tilde{E} (i \epsilon)^* \Pi_{\epsilon} a
\Pi_{\epsilon} \tilde{E}(i \epsilon)$}

 The following lemma will be used to reduce ergodicity properties
 and norm estimates of the complexified eigenfunctions
 $\phi_{\lambda}^{\C}$ to properties of the initial real
 eigenfunctions $\phi_{\lambda}.$ We recall that the $\tilde{E}, \tilde{\Pi}$
 notation refers to the $B^*_{\epsilon} M$ setting while $E, \Pi$ refers to the $M_{\epsilon}$ setting.

\begin{lem} \label{PSIDO} Let  $a \in S^0(T^*M-0)$.  Then for all $0 < \epsilon < \epsilon_0$,
we have: $$ \tilde{E} (i \epsilon)^* \tilde{\Pi}_{\epsilon} a
\tilde{\Pi}_{\epsilon} \tilde{E}(i \epsilon) \in \Psi^{-
\frac{m-1}{2}}(M),
$$
with principal symbol equal to $a(x, \xi) \; |\xi|_g^{- (
\frac{m-1}{2})}.$
\end{lem}

\begin{proof}

We observe that  $\tilde{E}(i \epsilon)^* \tilde{\Pi}_{\epsilon} a
\tilde{\Pi}_{\epsilon} \tilde{E}(i \epsilon)$ is a complex Fourier
integral operator on $L^2(M)$ associated to the canonical
relation,
$$\Gamma^* \circ \Delta_{\Sigma_{\epsilon}} \circ \Gamma =
 \{(y, \eta, x, \xi)\} \subset T^*M \times T^*M: \iota_{\epsilon}(x, \xi) =
\iota_{\epsilon} (y, \eta) \}.$$ Since $\iota_{\epsilon}$ is a
conic symplectic isomorphism, it follows that $\Gamma^* \circ
\Delta_{\epsilon} \circ \Gamma = \Delta_{T^*M \times T^*M}$, i.e.
that $\tilde{E}(i \epsilon)^* \tilde{\Pi}_{\epsilon} a
\tilde{\Pi}_{\epsilon} \tilde{E} (i \epsilon)$ is a
pseuododifferential operator. It follows from Theorem \ref{BDM}
that
$$\tilde{E}(i \epsilon)^* \tilde{\Pi}_{\epsilon} a \tilde{\Pi}_{\epsilon} \tilde{E}(i \epsilon): W^s(M) \to W^{s + \frac{m
- 1}{2}}(M),
$$
is a continuous linear map and hence that  the order is $-
\frac{m- 1}{2}.$

The  principal symbol of $ \tilde{E}(i \epsilon)^*
\tilde{\Pi}_{\epsilon} a \tilde{\Pi}_{\epsilon} \tilde{E}(i
\epsilon) $ equals $a(x, \xi)$ times the principal symbol of $
\tilde{E}(i \epsilon)^* \tilde{\Pi}_{\epsilon} \tilde{E}(i
\epsilon). $ Indeed, by the calculus and Egorov theorem for
complex Fourier integral operators, it equals  the the principal
symbol of $ \tilde{E}(i \epsilon)^* \tilde{\Pi}_{\epsilon}
\tilde{E}(i \epsilon)$ times the translate of $a$ under the
canonical relation underlying $\tilde{\Pi}_{\epsilon}
\tilde{E}(it)$. As discussed above, this relation is the
symplectic identification of $T^*M - 0 \equiv \Sigma_{\epsilon}.$

Thus, it suffices to show that  the  principal symbol of  $
\tilde{E}(i \epsilon)^*\tilde{\Pi}_{\epsilon}  \tilde{E}(i
\epsilon) $ equals $|\xi|_g^{-  \frac{m-1}{2}}.$  Since it equals
the principal symbol of $E(i t)^* \Pi_{\epsilon} E(it)$
transported under $\psi$,  we can (and will) do the computation in
the $M_{\epsilon}$ setting. The calculation could be done using
the calculus of complex Fourier integral operators, but that would
require a digression on symbols of \szego projectors and on the
composition of the symbols of the three factors. It seems quicker
to  calculate the symbol from scratch by applying stationary phase
for complex phase functions with only real critical points.

The calculation is well illustrated by   the simplest case of the
Euclidean wave kernel on $\R^m$. We then have:
$$E(i \tau)^* \Pi_{\tau} E(i \tau)(x, y) = \int_{\R^m} \int_{\R^m} \int_{\R^m \times S^{m-1}}
 e^{- i ( - i \tau)  |\xi_1|} e^{- i \langle \xi_1, x  - (x_1 - i p)
 \rangle}
  e^{i ( i \tau)  |\xi_2|} e^{i \langle \xi_2, x_1 + i p - y
\rangle} d\xi_1 d \xi_2 dp dx_1. $$ This is a complex Fourier
integral operator with phase
\begin{equation}\label{PHI0}  \Psi_0(x, y, \xi, \zeta, \tau) = - \tau (|\xi_1|
+ |\xi_2|) + i\langle \xi_2, x_1 + i p - y\rangle - i \langle
\xi_1, x - (x - i p) \rangle. \end{equation} The  $dx_1$ integral
produces $\delta(\xi_1 + \xi_2)$, and then  the $d\xi_1$ integral
gives
$$\begin{array}{lll} E(i \tau)^* \Pi_{\tau} E(i \tau)(x, y) & = & \int_{\R^m} \int_{\R^m} \int_{ S^{m-1}}
 e^{- i (t - i \tau)  |\xi|} e^{- i \langle \xi, x  - ( i p)
 \rangle}
  e^{i (t + i \tau)  |\xi|} e^{i \langle \xi,\ i p - y
\rangle} d\xi dp \\ & & \\
& = & \int_{\R^m}  \int_{ S^{m-1}}
 e^{- 2 \tau  |\xi|} e^{-  2 \langle \xi,   p
 \rangle}
 e^{- i \langle \xi, y
\rangle} d\xi dp \\ & &  \\
& \sim &  2 \tau^{\frac{m-1}{2}} \int_{\R^m}
|\xi|^{-\frac{m-1}{2}} \sinh 2 \tau  |\xi| e^{- 2 \tau |\xi|} e^{-
i \langle \xi, x - y \rangle} d\xi \end{array}$$ In the last step
we used the asymptotics of the inner integral
\begin{equation} \label{EUC} \begin{array}{lll} \int_{\tau S^{m-1}} e^{- 2
\langle \xi, p \rangle} d\mu(p) & = & \tau^{m-1} \int_{S^{m-1} }
e^{- 2 \tau |\xi| \langle \omega, e_1 \rangle} d\mu(\omega) \\ & &
\\ &\sim & 2 \tau^{\frac{m-1}{2}} |\xi|^{-\frac{m-1}{2}} \sinh 2
\tau  |\xi| e^{- 2 \tau |\xi|} .
\end{array} \end{equation} Thus, $E(i \tau)^* \Pi_{\tau} E(i \tau)$ is a
Fourier multiplier by a polyhomogeneous function with leading term
$|\xi|^{-\frac{m-1}{2}} $.

We now show that the same result holds on a general compact
analytic Riemannian manifold by reducing to  the Euclidean case.
Using the analytic continuation of the
 parametrix (\ref{CXPARAONE}), we have
\begin{equation} \label{CXPARATWO} \begin{array}{lll} E(i \tau)^* \Pi_{\tau} E(i \tau)(x, y)& =  &\int_{T^*_x M} \int_{T^*_y M} \int_{S^*M}
 e^{- \tau  (|\xi_1|_{g, x} +   |\xi_2|_{g, y} )}e^{- i \langle
\xi_1, \overline{\exp_x^{-1} (\zeta)} \rangle} e^{i \langle \xi_2,
\exp_y^{-1} (\zeta) \rangle} \\ & & \\ & \times & A^*(t, \zeta, x,
\xi_1) A(t, \zeta, y, \xi_2) d\xi_1 d \xi_2 d V(z) d \omega
\end{array}
\end{equation}

 The phase
\begin{equation} \label{PSI} \Psi = - \tau (|\xi_1|_{g,
\overline{\zeta}} + |\xi_2|_{g, y} ) - i \langle \xi_1,
\overline{\exp_x^{-1} (\zeta)} \rangle + i \langle \xi_2,
\exp_y^{-1} (\zeta) \rangle
\end{equation} is complex but has only real critical points given
by the non-degenerate critical manifold $C_{\psi} \simeq T^*M$:
\begin{equation} C_{\Psi} = \{(x, y, \xi_1, \xi_2, \zeta = z + i \tau \omega): y = x, \;\;\; \xi_1 = - \xi_2,\;  z = 0, \; \omega =
\frac{\xi_1}{|\xi_1|}\}. \end{equation} The  Lagrange immersion
$$\iota_{\Psi}: C_{\Psi} \to T^*M \times T^*M,\;\; \iota_{\Psi}(x, y, \xi_1, \xi_2, \zeta = z + i \tau
\omega)= (x, d_x \Psi, y, - d_y \Psi) $$ is then real valued and,
as mentioned above,  we may apply stationary phase for complex
phase functions with only real critical points (cf. \cite{Ho},
Vol. 1) to evaluate the principal symbol.

 We recall that
the  principal symbol  is the transport to $\Delta_{T^*M \times
T^*M}$ by $\iota_{\Psi}$ of the 1/2-density $\sigma_0
\sqrt{d_{C_{\Psi}}}$ where $\sigma_0$ is the principal term of the
amplitude restricted to $C_{\Psi}$, and where  $d_{C_{\Psi}}$ is
the Leray density on $C_{\Psi}$ induced by the map $i_{\Psi}$ and
the coordinate volume density $ d\xi_1 d \xi_2 d V(z) d \omega$.
Above, $\sigma_0 $ is the principal term of $A^*(t, \zeta, x,
\xi_1) A(t, \zeta, y, \xi_2) \equiv 1$ on $C_{\Psi}$. The Leray
density is given by
$$\begin{array}{l} d_{C_{\Psi}}  =  \left\| \frac{D(x, \xi_1,  \frac{\partial \Psi}{\partial \xi_1},
 \frac{\partial \Psi}{\partial \xi_2}, \frac{\partial \Psi}{\partial \zeta})}{D(x, \xi_1, y, \xi_2, \zeta)}\right\|^{-1} dx d\xi_1 =  \left\| \frac{D( \frac{\partial \Psi}{\partial \xi_1},
\frac{\partial \Psi}{\partial \xi_2}, \frac{\partial
\Psi}{\partial \zeta})}{D(y, \xi_2, \zeta)}\right\|^{-1} dx
d\xi_1,
\end{array}$$ where we regard $(x, \xi_1)$ as coordinates on
$C_{\Psi}$. We write $\zeta = x_1 + i p$ in Riemannian normal
coordinates based at $x$. Expanding first along the $D_y$ rows and
then along the $\Psi'_{x_1}$ columns, we obtain
\begin{equation} \label{END} \begin{array}{l} d_{C_{\Psi}}   =  \left\| \frac{D(
\frac{\partial \Psi}{\partial p)}}{D(p)}\right\|^{-1} dx d\xi_1.
\end{array} \end{equation}

Indeed, simple calculations show that the $D_y$ derivatives only
act non-trivially on $\Psi_{\xi_2}'$, and that $\Psi'_{x_1}$ only
has non-trivial derivatives under $D_{\xi_2}, D_{x_1}.$ To
eliminate the $D_y$ rows and the $\Psi'_{x_1}$ column, it suffices
to show that $\Psi''_{x_1 x_1} = 0$ and that $\Psi''_{y \xi_2} =
Id = \Psi''_{x_1 \xi_2}$ on the critical set. In the calculation
of $\Psi''_{x_1 x_1}$ on the critical set we may put $\xi_1 = -
\xi_2 = \xi, x = y, p = \frac{\xi_1}{|\xi_1|}$, and then the
calculation reduces to $D^2_{x_1} (\langle \xi, \exp_x^{-1} (x_1 -
i p) - \exp_x^{-1}(x_1 + i p)\rangle)|_{x_1 = 0} = 0,$ which holds
since $\exp_x^{-1}$ is the identity map when $\zeta = x + i p$ is
expressed in Riemannian normal coordinates at $x$. In calculating
$\Psi''_{y \xi_2} = D^2_{y \xi_2} (|\xi_2|_{g_y} + i \langle
\xi_2, \exp_y^{-1}(\zeta) \rangle) |_{x = y, \xi_2 = - \xi_1}$,
the first term has a zero Hessian since all metric coefficients
vanish in normal coordinates. The second is easily seen to be a
constant multiple of the identity operator (the multiple is the
same as in the Euclidean case). Finally, $\Psi''_{x_1 \xi_2} = i
D_{x_1} \exp_y^{-1} (x_1 + i p)|_{x_1 = 0, y = x} = C Id$. Taking
the determinant of the result gives (\ref{END}). Finally, we
calculate this last determinant after restricting to the critical
set, obtaining a constant multiple of
$$\det D^2_{pp} \langle \xi, \exp_x^{-1}(ip) \rangle |_{p =
\frac{\xi}{|\xi|}} = C \det D^2_{pp} \langle \xi, p \rangle |_{p =
\frac{\xi}{|\xi|}}. $$ The last determinant is precisely the one
that arises as the Hessian determinant in  the stationary phase
formula in (\ref{EUC}). The derivatives are taken on $S^{m-1}$,
hence we obtain $|\xi|^{m-1}$ times the normalized determinant,
which is invariant under rotations and hence constant. Raising to
the power $-\frac{1}{2}$ completes the calculation.

\end{proof}

\section{Proof of Lemma  \ref{ERGOCOR} and  Lemma
\ref{NORM}}

We now use Lemma \ref{PSIDO} to reduce the quantum ergodicity and
norm properties of the complexified eigenfunctions to properties
of the original real eigenfunctions.

\subsection{Proof of Lemma \ref{ERGOCOR}}

  We begin by proving a
 weak limit formula for the CR holomorphic functions $u_{\lambda}^{\epsilon}$  defined in
 (\ref{LITTLEU}) for fixed $\epsilon$. For notational simplicity,
 we drop the tilde notation although we work in the
 $B^*_{\epsilon} M$ setting.

\begin{lem} \label{ERGO} Assume that $\{\phi_{\lambda}\}$ is a quantum ergodic sequence.
Then for each $0 < \epsilon < \epsilon_0$,
$$|u_{\lambda}^{\epsilon}|^2 \to
\frac{1}{\mu_{\epsilon}(\partial B^*_{\epsilon} M)},\;\;
\mbox{weakly in}\;\; L^1(\partial B^*_{\epsilon}M,
d\mu_{\epsilon}).
$$
That is, for any   $a \in C(\partial B^*_{\epsilon}M)$,
$$\int_{\partial B^*_{\epsilon} M} a(x, \xi) |u_{\lambda}^{\epsilon}((x, \xi)|^2
d\mu_{\epsilon} \to \frac{1}{\mu_{\epsilon}(\partial
B^*_{\epsilon} M)} \int_{\partial B^*_{\epsilon}M} a(x, \xi)
d\mu_{\epsilon}.$$

\end{lem}

\begin{proof}

It suffices to consider $a \in C^{\infty} (\partial
B^*_{\epsilon}M)$.  We then consider the Toeplitz operator
$\Pi_{\epsilon} a \Pi_{\epsilon}$ on $\ocal^0(\partial
B^*_{\epsilon} M)$. We have,
\begin{equation}\label{EQUIV} \begin{array}{lll}  \langle \Pi_{\epsilon} a \Pi_{\epsilon}
u_j^{\epsilon}, u_j^{\epsilon} \rangle & =  & e^{2 \epsilon
\lambda_j}\;  ||\phi_{\lambda}^{\C}||_{L^2(\partial B^*_{\epsilon}
M)}^{-2} \langle \Pi_{\epsilon} a \Pi_{\epsilon} E(i \epsilon)
 \phi_j, E(i \epsilon) \phi_j
\rangle_{L^2(\partial B^*_{\epsilon} M)}\\  & & \\
& = &  e^{2 \epsilon \lambda_j}\;
||\phi_{\lambda}^{\C}||_{L^2(\partial B^*_{\epsilon} M)}^{-2}
\langle E(i \epsilon)^* \Pi_{\epsilon} a \Pi_{\epsilon} E(i
\epsilon) \phi_j, \phi_j \rangle_{L^2(M)}.
\end{array}
\end{equation}
 By
Lemma \ref{PSIDO}, $ E(i \epsilon)^* \Pi_{\epsilon} a
\Pi_{\epsilon} E(i \epsilon) $ is a pseudodifferential operator on
$M$ of order $- \frac{m-1}{2}$ with principal symbol $\tilde{a}
|\xi|_g^{ - \frac{m-1}{2}}$, where $\tilde{a}$ is the (degree $0$)
homogeneous extension of $a$ to $T^*M - 0$.  The normalizing
factor $ e^{2 \epsilon \lambda_j}\;
||\phi_{\lambda}^{\C}||_{L^2(\partial B^*_{\epsilon} M)}^{-2}$ has
the same form with $a = 1$. Hence, the expression on the right
side of  (\ref{EQUIV}) may be written as
\begin{equation} \frac{\langle E(i \epsilon)^* \Pi_{\epsilon} a \Pi_{\epsilon} E(i
\epsilon) \phi_j, \phi_j \rangle_{L^2(M)}}{\langle E(i \epsilon)^*
\Pi_{\epsilon} E(i \epsilon) \phi_j, \phi_j \rangle_{L^2(M)}}.
\end{equation}

  By the standard quantum ergodicity result on compact
Riemannian manifolds with ergodic geodesic flow (see \cite{Shn,Z2,
Z, CV} for proofs and references)  we have
\begin{equation} \label{QE} \frac{\langle E(i \epsilon)^* \Pi_{\epsilon} a \Pi_{\epsilon} E(i
\epsilon) \phi_j, \phi_j \rangle_{L^2(M)}}{\langle E(i \epsilon)^*
\Pi_{\epsilon} E(i \epsilon) \phi_j, \phi_j \rangle_{L^2(M)}} \to
\frac{1}{\mu_{\epsilon}(\partial B^*_{\epsilon} M)} \int_{\partial
B^*_{\epsilon} M} a d \mu_{\epsilon}.
\end{equation}
More precisely, the numerator is asymptotic to the right side
times $\lambda^{- \frac{m-1}{2}}$, while the denominator has the
same asymptotics when $a$ is replaced by $1$. We also use that
$\frac{1}{\mu_{\epsilon}(\partial B^*_{\epsilon} M)}
\int_{\partial B^*_{\epsilon} M} a d \mu_{\epsilon}$ equals the
analogous average of $\tilde{a}$ over $\partial B_1$ (see the
discussion around (\ref{LIOUVSYM})). Taking the ratio produces
(\ref{QE}).

Combining (\ref{EQUIV}), (\ref{QE}) and the fact that
$$ \langle \Pi_{\epsilon} a \Pi_{\epsilon}
u_j^{\epsilon}, u_j^{\epsilon} \rangle = \int_{\partial
B^*_{\epsilon} M} a |u_j^{\epsilon}|^2 d \mu_{\epsilon}$$
completes the proof of the lemma.

\end{proof}

We now complete the proof of Lemma \ref{ERGOCOR}, i.e. we prove
that
\begin{equation} \label{WEAKLIMA} \int_{B^*_{\epsilon} M} a
|U_{\lambda}|^2 \omega^m \to \frac{1}{\mu_1(S^* M)}
\int_{B^*_{\epsilon} M} a |\xi|_g^{-m + 1} \omega^m
\end{equation} for any $a \in C(B^*_{\epsilon} M).$
It is only necessary to relate the Liouville measures $d\mu_r$
(\ref{LIOUVILLE}) to the symplectic volume measure. One may write
$d\mu_r = \frac{d}{dt} |_{t = r} \chi_t \omega^m$, where $\chi_t$
is the characteristic function of $B_t^*M = \{ |\xi|_g \leq t\}$.
By homogeneity of $|\xi|_g$, $\mu_{r}(\partial B^*_{r} M) =
r^{m-1} \mu_1(\partial B^*_1 M) $. If  $a \in C(B^*_{\epsilon})$,
then $\int_{B_{\epsilon}^*M} a \omega^m = \int_0^{\epsilon}
\{\int_{\partial B^*_{r} M} a
 d \mu_{r}\} dr.  $ By Lemma \ref{ERGO}, we  have

\begin{equation} \label{LIOUVSYM} \begin{array}{lll} \int_{B_{\epsilon}^*M} a |U_{\lambda}|^2
\omega^m = \int_0^{\epsilon}\{ \int_{\partial B^*_{r} M} a
|u_{\lambda}^r|^2  d \mu_{r} \}  d r &\to & \int_0^{\epsilon}
\{\frac{1}{\mu_{r} (\partial B^*_{r})} \int_{\partial B^*_{r} M} a
d \mu_r \}
dr \\ & & \\
& = &  \frac{1}{\mu_1(\partial
 B_1^* M)}\int_{B^*_{\epsilon} M} a r^{-m + 1} \omega^m, \\ & & \\
& \implies & w^*-\lim_{\lambda \to \infty} |U_{\lambda}|^2 =
 \frac{1}{\mu_1(\partial B_1^* M)} |\xi|_g^{-m + 1}.
\end{array}
\end{equation}

\subsection{Proof of Lemma \ref{NORM}}

We actually prove the stronger result that
$$\frac{1}{\lambda} \log \rho_{\lambda}(x, \xi)  \to |\xi|_{g_x},\;\;\;\mbox{uniformly in }\;\;B^*_{\epsilon} M\;\; \mbox{as}\;\; \lambda \to
\infty.
$$
We state the weaker form because that is what we need in the proof
of Theorem \ref{ZERO}.

Again we drop the tilde notation for simplicity.

\begin{proof}

Again using  $E(i \epsilon) \phi_{\lambda} = e^{-\lambda \epsilon}
\phi_{\lambda}^{\C}$, we have:
 \begin{equation} \label{ME} \begin{array}{lll} \rho^2_{\lambda}(x, \xi) &
 =&
 \langle \Pi_{\epsilon} \phi_{\lambda}^{\C}, \Pi_{\epsilon} \phi_{\lambda}^{\C} \rangle_{L^2(\partial B^*_{\epsilon} M)} \;\;(\epsilon = |\xi|_{g_x})\\ & & \\
 & = &
 e^{  2 \lambda
\epsilon}
 \langle  \Pi_{\epsilon} E(i \epsilon) \phi_{\lambda}, \Pi_{\epsilon} E(i \epsilon)
 \phi_{\lambda}
\rangle_{L^2(\partial B^*_{\epsilon} M)} \\ & & \\ & = &   e^{ 2
\lambda \epsilon} \langle E(i \epsilon)^* \Pi_{\epsilon} E(i
\epsilon) \phi_{\lambda},  \phi_{\lambda} \rangle. \end{array}
\end{equation}
Hence,
\begin{equation} \label{LOGAR} \frac{2}{\lambda} \log
\rho_{\lambda} (x, \xi)= 2 |\xi|_{g_x} + \frac{1}{\lambda} \log
\langle E(i \epsilon)^* \Pi_{\epsilon} E(i \epsilon)
\phi_{\lambda}, \phi_{\lambda} \rangle. \end{equation}

To complete the proof, we observe that
 \begin{equation} \label{LOGBOUND} \frac{1}{\lambda} \log \langle E(i \epsilon)^* \Pi_{\epsilon} E(i \epsilon)
\phi_{\lambda}, \phi_{\lambda} \rangle \leq C\; \frac{\log
\lambda}{\lambda}, \;\; \mbox{uniformly in } \; \epsilon,
\end{equation}
where $C$ is a constant independent of $(\epsilon, \lambda)$.
Indeed, by Lemma \ref{PSIDO}, $ E(i \epsilon)^* \Pi_{\epsilon} E(i
\epsilon)$ is a pseudodifferential operator of order $- \frac{m-
 1}{2}$ with principal symbol $|\xi|^{-\frac{m-1}{2}}.$ To obtain
 a uniform bound on $\langle E(i \epsilon)^* \Pi_{\epsilon} E(i \epsilon)
\phi_{\lambda}, \phi_{\lambda} \rangle$ in $\epsilon$,  any  of
the standard bounds for the norm of a pseudo-differential operator
in terms of derivatives of the complete symbol would suffice.

To take one such bound with a convenient reference, we write
$\langle E(i \epsilon)^* \Pi_{\epsilon} E(i \epsilon)
\phi_{\lambda}, \phi_{\lambda} \rangle$ as $(1 +
\lambda^2)^{\frac{m+1}{2}}
 \langle (I - \Delta )^{-\frac{m+1}{2}} E(i \epsilon)^* \Pi_{\epsilon} E(i \epsilon)
\phi_{\lambda}, \phi_{\lambda} \rangle$. Put $A_{\epsilon}:=  (I -
\Delta )^{-\frac{m+1}{2}} E(i \epsilon)^* \Pi_{\epsilon} E(i
\epsilon)$. Since $A_{\epsilon} \in \Psi^{-(m-1)}$, we may apply
the Schur-Young bound of \cite{H} (Vol. III, Theorem 18.1.11) to
obtain
\begin{equation}\label{SCHUR}  \langle A_{\epsilon} \phi_{\lambda},
\phi_{\lambda} \rangle  \leq ||A_{\epsilon}||_{L^2 \to L^2} \leq
C_m \left(\sup_x \int_{T^*_x M} |a_{\epsilon}(x, \xi)| d \xi
\right),
\end{equation} where $a_{\epsilon}$ is the complete symbol
of $A_{\epsilon}$ relative to some choice of quantization $a(x,
D)$ of symbols. The complete symbol of $A_{\epsilon}$ may be
obtained by applying $ \langle (I - \Delta )^{-\frac{m+1}{2}}$ to
the representation in (\ref{CXPARATWO}). It is clear that the
complete symbol is smooth in the parameter $\epsilon$, hence the
right side of (\ref{SCHUR}) has uniform bound in $\epsilon$,
proving (\ref{SCHUR}) and therefore the lemma.

\end{proof}

\section{\label{PROOF}  Proof of Lemma \ref{ZEROWEAK} and Theorem \ref{ZERO}}

The remaining  step in the proof of Theorem \ref{ZERO} is  the
proof of Lemma \ref{ZEROWEAK}.

\subsection{Proof of Lemma \ref{ZEROWEAK}}

\begin{proof} The proof is similar to that of Lemma 1.4 of
\cite{SZ}. We wish to prove that  $$\psi_j := \frac{1}{\lambda_j}
\log |U_j|^2 \to 0 \;\; \mbox{in}\; L^1(B^*_{\epsilon} M).
$$
We argue by contradiction. If the conclusion is not true, then
there exists a subsquence $\psi_{j_k}$ satisfying
$||\psi_{j_k}||_{L^1(B^*_{\epsilon} M)} \geq \delta > 0.$

To obtain a contradiction, we first observe that $\psi_j$ is
quasi-plurisubharmonic
 (QPSH)
 on $B^*_{\epsilon} M$, i.e. may be locally written as the sum of a plurisubharmonic
 function $v_j$ and a smooth function $\rho_j$;  equivalently  $i \ddbar \psi_j$ is locally bounded below by a negative
 smooth $(1,1)$ form. Indeed we put
$$v_j := \frac{1}{\lambda_j} \log |\phi_j^{\C}|^2,\;\; \rho_j
:= - \rho_{\lambda_j}.$$  We   use the following fact about
subharmonic functions (see \cite[Theorem~4.1.9]{Ho}):
\medskip

{\it Let $\{v_j\}$ be a sequence of subharmonic functions in an
open set $X \subset \R^m$ which have a uniform upper bound on any
compact set. Then either $v_j \to -\infty$ uniformly on every
compact set, or else there exists a subsequence $v_{j_k}$ which is
convergent in $L^1_{loc}(X)$. }

 Since the proof is
local, it also holds for open sets in manifolds, and in particular
for $X = B^*_{\epsilon} M$.

 We now verify that  the hypotheses are satisfied in our example:

\begin{itemize} \item (i)  the functions $v_j$ are
uniformly bounded above on $B^*_{\epsilon} M$;
\item (ii)  $\limsup_{j \rightarrow \infty}  v_j \leq 2 |\xi|_g $.
\end{itemize}

It suffices to prove these statements on each surface $\partial
B^*_{\epsilon} M$ with uniform constants independent of
$\epsilon$. On the surface $\partial B^*_{\epsilon} M$, $U_j =
u^{\epsilon}_j$. By the Sobolev inequality in
$\ocal^{\frac{m-1}{4}}(\partial B^*_{\epsilon} M)$, we have
$$\begin{array}{lll} \sup_{(x, \xi) \in \partial B^*_{\epsilon} M)} |u_j^{\epsilon} (x, \xi)| & \leq &
\lambda_j^m ||u_j^{\epsilon} (x, \xi)||_{L^2(\partial
B^*_{\epsilon} M)} \\ & & \\
& \leq & \lambda_j^m.
\end{array}$$
 Taking
the logarithm, dividing by $\lambda_j$, and combining with the
limit formula of Lemma \ref{NORM} proves (i) - (ii).

We now settle the dichotomy above by proving that the sequence
$\{v_j\}$ does not tend uniformly to $-\infty$ on compact sets.
That would imply that $\psi_j \to - \infty$ uniformly on the
spheres $\partial B^*M_{\epsilon}$ for each $\epsilon <
\epsilon_0$. Hence, for each $\epsilon$, there  would exist $K
> 0$ such that for $k \geq K$,
\begin{equation}\frac{1}{\lambda_{j_k}} \log | u_{j_k}^{\epsilon} (z)| \leq -
1.\label{firstposs}\end{equation} However, (\ref{firstposs})
implies that
$$ | u_{j_k}(z)| \leq e^{- 2 \lambda_{j_k}}\;\;\;\;\forall z \in
\partial B^*_{\epsilon} M\;, $$ which is inconsistent with the hypothesis that $|
u_{j_k}^{\epsilon} (z)| \rightarrow 1$ in $\dcal'(\partial
B^*_{\epsilon} M)$.

Therefore, the second half of the dichotomy holds, i.e.
 there must exist a subsequence, which we continue to denote by
$\{v_{j_k}\}$, which converges in $L^1(B^*_{\epsilon_0})$ to some
$v \in L^1(B^*_{\epsilon_0}).$ By passing if necessary to a
further subsequence, we may assume that $\{v_{j_k} \}$ converges
pointwise almost everywhere to $v$ in $B^*_{\epsilon_0}$. Then,
$$v (z) = \limsup_{k \rightarrow \infty} v_{j_k}\leq 2 |\xi|_g\;\;\;\;\;\;{\rm (a.e)}\;.$$ Now let
$$v^*(z):= \limsup_{w \rightarrow z} v(w) \leq 0 $$ be the
upper-semicontinuous regularization of $v$. Then $v^*$ is
plurisubharmonic on $B^*_{\epsilon} M$ and $v^* = v$ almost
everywhere.

Put $\psi^* : = v^* - 2 |\xi|_g.$ Then $\psi^* \leq 0$, and the
assumption $||\psi_{j_k}||_{L^1(B^*_{\epsilon} M)} \geq \delta >
0$ implies that
$$U_{\delta} : = \{\zeta \in B^*_{\epsilon_0} M: \; \psi^* (\zeta)\; <  -
\delta/2\} $$ has positive volume. Since $\psi_{j_k} \to \psi^*$
in $L^1(U_{\delta})$, one has by \cite{Ho} Theorem 4.1.9 (b)  that
\begin{equation} \limsup_{k \to \infty} \psi_{j_k} |_{U_{\delta}}
\leq \psi^*|_{U_{\delta}} < - \delta/2.\end{equation} Hence, there
exists a positive integer $K$ such that $\psi_{j_k}(\zeta) \leq
-\delta/2$ for $\zeta \in U_{\delta},\ k\geq K$; i.e.,
\begin{equation} |\psi_{j_k}(\zeta)|\leq e^{- \delta
\lambda_{j_k}},\;\;\;\;\; \zeta\in U_{\delta}, \;\;k \geq K.
\end{equation} This again contradicts the weak convergence to $1$.

 Therefore $||\psi_{j_k}||_{L^1(B^*_{\epsilon} M)} \geq \delta >
0$ leads to a contradiction, and the Lemma is proved.

\end{proof}

To complete the proof of Theorem \ref{ZERO} it suffices to combine
the results that
$$\frac{i}{2 \pi \lambda} \ddbar \log |U_{\lambda}|^2 = \frac{1}{\lambda}\; [Z_{\lambda}] - \frac{i}{2 \pi \lambda}
\ddbar \log ||\phi_{\lambda}^{\epsilon}||^2_{L^2(\partial
M_{\epsilon})} \to 0 \;\;\;(\mbox{weakly}),  $$ and that  (by
Lemma \ref{NORM}) the second term tends to $ \frac{i}{ \pi} \ddbar
|\xi|_g.$ \qed

 \subsection{\label{FINAL}
 Final remarks}

 \medskip

\noindent{\bf (i)} We check the numerical details in the case of
the circle.

 The zeros of $\sin 2 \pi k z$ in the cylinder
 $\C/\Z$ all lie on the real axis at the points $z = \frac{n}{2
 k}.$ Thus, there are $2 k$ real zeros, and  the Poincar\'e-Lelong
 formula gives
 $$\begin{array}{lll} \lim_{k \to \infty} \frac{i}{2 \pi k} \ddbar \log |\sin 2 \pi k|^2
& = &\lim_{k \to \infty} \frac{1}{k} \sum_{n = 1}^{2 k} \delta_{\frac{n}{2 k}} \\ & & \\
& = &
 \frac{ 1}{\pi}  \delta_0(\xi) dx \wedge d \xi. \end{array}$$
 On the other hand,
 $$\begin{array}{lll} \frac{i}{ \pi} \ddbar |\xi| & = & \frac{i}{ \pi}
 \frac{d^2}{4 d\xi^2} |\xi| \;\; \frac{2}{i} dx \wedge d\xi \\ & & \\
 & = & \frac{i}{ \pi}\;
 \frac{1}{2} \;   \delta_0(\xi) \;\; \frac{2}{i} dx \wedge d\xi, \end{array}$$
 matching the other expression.

 As mentioned in the introduction,  the complex eigenfunctions $e^{2 \pi i k x}$ have no
 complex zeros, hence Theorem \ref{ZERO} is false for them. The
  reason is that they are not quantum ergodic but rather
 localize on just one of the two components of the unit tangent
 bundle (the one with the same sign as $k$). Running through the
 previous calculation shows that the limit zero current for these eigenfunctions is $\frac{i}{ \pi} \ddbar
 \xi = 0$ rather than $\frac{i}{ \pi} \ddbar |\xi|$.
 \medskip

 \noindent{\bf (ii)} One can obtain other formulae for the
 distribution of zeros in the ergodic case using the fact that the
 maps $t + \sqrt{-1} s \to \exp_{\gamma(t)} s
\dot{\gamma}(t) $ are holomorphic curves  $\gamma_{\C}(t +
\sqrt{-1} s)$ relative to the adapted complex structure for each
geodesic $\gamma$. If one pulls back complexified eigenfunctions
under $\gamma_{\C}$, then one obtains a holomorphic function in a
strip around the real-axis. Its complex zeros are  discrete and
correspond to the intersection points $Z_{\phi_{\lambda}^{\C}}
\cap \gamma_{\C}$.  Its  real zeros are the intersection points of
the real geodesic $\gamma$ with the nodal hypersurface. In
connection with Conjecture \ref{REALZ}, it is natural to
conjecture that these intersection points become uniformly
distributed on $(M, g)$ when $\gamma$ is a uniformly distributed
geodesic.
\medskip

\noindent{\bf (iii)} We can give  a simpler form to $ \frac{i}{
\pi} \ddbar
 |\xi|_g$ in dimensions $m \geq 2.$
Let $\rho(x, \xi) = |\xi|_g^2$. We  note (with \cite{GS2}) that
 $$\ddbar f(\rho) = f'(\rho) \ddbar \rho + f''(\rho)
 \overline{\partial} \rho \wedge \partial \rho, $$
 and that
 $$\overline{\partial} \rho \wedge \partial \rho = i d \rho \wedge
 \alpha,\;\;\; \ddbar \rho = - i \omega_g. $$
 It follows that in dimensions $m \geq 2,$ we have
 \begin{equation} \frac{i}{\pi} \ddbar \sqrt{\rho} = \frac{1}{2 \pi \rho^{1/2}}
 \omega_g + \frac{ d \rho \wedge \alpha}{4 \pi  \rho^{3/2}}.
 \end{equation}

\end{document}